\documentclass[10pt]{article}
\usepackage{amsmath}
\usepackage{amsfonts}
\usepackage{amscd}
\usepackage{graphicx}
\usepackage{amssymb}
\usepackage{latexsym}
\usepackage{float}
\usepackage{color}

\newcommand{\R}{\mathbb{R}}

\def \r {\color{red}}
\def \d {\displaystyle}
\newtheorem{thm}{Theorem}
\newtheorem{lem}{Lemma}

\author{\bf Berge Tsanou$^{1,5}$, Samuel Bowong$^{2,5,\dag}$, Jean Jules Tewa$^{3,5}$, Gauthier Sallet$^4$ \\
\\
$^1$Department of Mathematics and Computer Science, Faculty of Science,\\
 University of  Dschang, P.O. Box 67 Dschang, Cameroon\\
 E-mail: bergetsanou@yahoo.fr\\
\\
$^2$Laboratory of Applied Mathematics, \\
Department of Mathematics and Computer Science, Faculty of Science, \\
University of Douala, P.O. Box 24157  Douala, Cameroon\\
E-mail: sbowong@gmail.com \\
\\
$^3$Department of Mathematics and Physics,\\
  National Advanced School of Engineering (Polytechnic), \\
  University of Yaounde I, P.O. Box 8390 Yaounde, Cameroon\\
  E-mail: tewajules@yahoo.fr \\
 \\
$^4$ MASAIE project-team, INRIA Grand-Est, and University Paul Verlaine-Metz,\\
LMAM-CNRS,UMR 7122 ISGMP Bat. A, Ile du Saulcy, 57045 Metz Cedex 01, France\\
\\
 $^5$ UMI 209 IRD/UPMC UMMISCO, Bondy, France and  \\
 Project-Team GRIMCAPE, LIRIMA, University of Yaounde I, Cameroon\\
 \\
$^\dag$ Corresponding author: Tel. +(237) 99-96-41-64, Fax. +(237) 22-31-02-90\\
E-mail: sbowong@gmail.com}%
\date{\empty}
 \title{\bf Analysis of the spread of tuberculosis in heterogeneous
complex metapopulations}
 
\begin{document}
\baselineskip7mm
\maketitle%
\begin{abstract} This paper describes and  analyzes  the spatial spread of tuberculosis (TB) on complex metapopulation, that is, networks of populations connected by migratory flows whose configurations are described in terms of connectivity distribution of nodes (patches) and the conditional probabilities of connections among classes of nodes sharing the same degree. The migration and transmission processes occur simultaneously.
For uncorrelated networks under the assumption of standard incidence transmission, we compute the disease-free equilibrium and the basic reproduction number, and show that the disease-free equilibrium is locally asymptotically stable. Moreover, for uncorrelated networks and under assumption of simple mass action  transmission, we give a necessary and sufficient conditions for the instability of the disease-free equilibrium.
 The existence of endemic equilibria is also discussed. Finally,  the prevalence of the TB infection across the metapopulation as a function
of the path connectivity is studied using numerical simulations.
 \end{abstract}


\section{Introduction}
\noindent

 Despite significant
advances in medical science, infectious diseases continue to impact human populations
in many parts of the world.  Tuberculosis (abbreviated as TB for tubercle bacillus) is a common deadly infectious disease caused mainly by {\it Mycobacterium
tuberculosis}. It basically attacks the lungs (pulmonary TB), but can also affect
the central nervous system, circulatory system, the genital-urinary system, bones,
joints and even the skin. Tuberculosis can spread through cough, sneeze, speak,
kiss or spit from active pulmonary TB persons. It can also spread through use of
an infected person's unsterilized eating utensils and in rare cases a pregnant woman with active TB can infect her foetus (vertical transmission) [1,2]. Transmission can
only occur from people with active TB but not latent TB. This transmission from
one person to another depends upon the number of infectious droplets expelled by a
carrier, the effectiveness of ventilation, duration of the exposure and virulence of the {\it MTB} strain. The chain of transmission can therefore be broken by isolating patients with active disease and starting effective anti-tuberculosis therapy [1-5]. At present, about $95\%$ of the estimated 8 million new cases of  TB occurring each year are in developing countries,
where $80\%$ occur among people between the ages of 15-59 years [1]. In sub-Saharan Africa, TB is the leading cause of mortality and in developing countries, it accounts for an estimated 2 million deaths which accounts for a quarter of avoidable adult deaths [1]. It is known that factors such as endogenous reactivation, emergence of multi-drug resistant TB, and increase in HIV incidence in the recent years call for improved control strategies for TB.  A full understanding of
the effectiveness of treatment and control strategies within different regions of the world is still needed. It is worth emphasizing that mathematical analysis of biomedical and disease transmission models can contribute
to the understanding of the mechanisms of those processes and to design potential therapies (see [6-9] and references therein). A number of theoretical studies have been carried out on the mathematical modeling
of TB transmission dynamics [3-9,38,39].

However, the analysis of the spread of infectious diseases on complex networks has become a central issue in modern epidemiology [10] and, indeed, it was one of the main motivations for the development of percolation theory [11]. While the initial approach was focussed on local contact networks [12-16], that is, social networks within single populations (cities, urban areas), a new approach has been recently introduced  for dealing with the spread of diseases in ensembles of (local) populations with a complex spatial arrangement and connected by the migrations.
Such sets of connected populations living in a patchy environment are called
metapopulations in ecology, and their study began in 1967 with the theory of island biogeography [17].

Unfortunately,  when considering dispersal models, there is an approach based on the
metapopulation concept. The population is subdivided into a number of discrete patches
which are supposed to be well mixed. Then, in each patch the population is
subdivided into compartments corresponding to different epidemic status. This
leads to a multi-patch, multi-compartment system. At this point two formulations
are possible.

  The first one assumes that an infective in one patch can infect susceptible
individuals in another patch. This assumption gives rise to a family of models
which have been well studied [18,19].  This formulation assumes that
there is a spatial coupling between patches, but that individuals (vectors or
hosts) do not migrate between patches. They make short `visits' from their
home patches to other patches.

The second one considers migration of individuals between patches. The infection
does not take place during the migration process. The situation is that
of a directed graph, where the vertices represent the patches and the arcs represent
the links between patches. Recently, there has been increased interest
in these deterministic metapopulation disease models. For instance, in some recent models of epidemic spreading, the location of the patches in space is treated  explicitly (without taking into account  the number of connections $ k $  (degree) that any given patch in the network may have) thanks to the increasing of computational power (see for instance [20, 21]). In Refs.~[16, 22, 23], however, an alternative approach based on the formalism used in the statistical mechanics of complex networks is presented. Under this approach, the structure of the spatial network of patches (nodes) is encapsulated by means of the connectivity (degree) distribution $p(k)$  defined as the probability that a randomly chosen patch has connectivity $k$. In  contrast, in [24, 25], the authors consider reaction diffusion processes  to take place simultaneously, which turns out to be correct assumption for a suitable continuous-time formulation of metapopulation models for the spread of infectious diseases.

In this paper, we consider the spread of TB on complex metapopulations, that is, networks of populations connected by migratory flows whose configurations are described in terms of the conditional probabilities of connections among classes of nodes sharing the same degree. For uncorrelated networks under the assumption of standard incidence [37] (or frequency-dependent) transmission, we compute the disease-free equilibrium and the basic reproduction number and show that the disease-free equilibrium is locally asymptotically stable. Moreover, for uncorrelated networks and under assumption of simple mass action [37] (or density-dependent) transmission, we give a necessary and sufficient conditions for the instability of the disease-free equilibrium.  We find that there exists a more precise bound of the largest eigenvalue of the Jacobian matrix of the system around the disease-free equilibrium. This condition says that, for fixed values of the migration rates of latently-infected and infectious individuals, a high enough density of individuals and/or large enough maximum connectivity in the metapopulation guarantee the instability of the disease-free equilibrium and, hence, TB spread. In the limit of infinite networks with bounded average degree, this condition implies the existence of a TB threshold for any distribution with large value. The existence of endemic equilibria is also discussed.   Additionally, through numerical simulations, the forecasted prevalence of the infection is not constant but increases with the patch connectivity. Interestingly, close the epidemic threshold, there are always patches with low connectivities where TB is not able to progress unless infectious individuals arrive from (crowded) patches with higher connectivities.  Comparing to existing results in the literature, our work treats a specific disease which is not the case in  Refs. [24, 25, 26]. We point out that in Refs. [24, 25, 26],  the authors have  neglected some important epidemiological features of the propagation of a disease  such as births, natural mortality,  mortality due to the disease, natural recovery and  the basic models studied are of dimension 2 which are very simple. In addition, the authors have supposed that the total population is constant which is not always the case. Our basic model is of dimension 4 and  incorporates the essential biological and epidemiological features of TB such as birth, mortality due to the disease, slow and fast progression, effective chemoprophylaxis of latently-infected individuals,  natural recovery and  treatment of infectious, relapse from the disease  and re-infection after recovery.   Also in our model, the total population is not constant. It is our view fact that this study represents the first work that provides an in-depth the spread  of TB on complex metapopulation using a degree of distribution and conditional probabilities.
\section{ A TB metapopulation model}
\noindent

\subsection{The model}
\noindent

We consider the spread of TB in heterogeneous  metapopulations.
The model consists of $n$ patches representing $n$ different degree of connectivities. We assume that the architecture of the network of patches (nodes) where local populations live is mathematically encoded by means of the connectivity (degree) distribution $p(k)$, defined as the probability that a randomly chosen patch has degree $k$. At any given time, in each patch, an individual is in one of the following states: susceptible, latently infected (exposed to TB but not infectious), infectious (has active TB) and recovered. These states are average number (density) of  $\rho_{S,k}$, $\rho_{E,k}$ , $\rho_{I,k}$ and  $\rho_{R,k}$   in the patches of connectivity $k$, respectively. The total variable population size at time $t$ is given by,
\begin{equation}
 \label{pop}
\rho_k(t) = \rho_{S,k}(t)+\rho_{E,k}(t)+\rho_{I,k}(t)  + \rho_{R,k}(t).
\end{equation}
 It is assumed that births
are recruited into the population at per capita rate $\Lambda$. The transmission of {\it Mycobacterium tuberculosis} occurs following adequate contacts between a susceptible and infectious  in each sub-population. The rate at which susceptible are infected is $\beta\displaystyle\frac{\rho_{I,k}\rho_{S,k}}{\rho_k}$ for standard incidence (or frequency-dependent) transmission and $\beta\,\rho_{I,k}\rho_{S,k}$  for simple mass action (or density-dependent) transmission,
where $\beta$ is  the effective contact rate of  infectious that is sufficient to transmit
infection to susceptible (it also denotes how contagious of the disease is). On adequate contacts with  active  individuals, a susceptible individual becomes infected but not yet infectious.
A fraction $q$  of newly infected individuals is assumed to undergo a fast progression
directly to the  infectious class, while the remainder is latently infected and enter the latent class.
Latently infected individuals are assumed to acquire some immunity as a result of the infection, which reduces the risk of
subsequent infection but does not fully prevent it.  We assume that
chemoprophylaxis of latently infected individuals reduces
their reactivation at a  rate $\theta$ and that the initiation of
therapeutics immediately remove individuals from active status and place them into a latent state.  This last assumption is realistic. Indeed, the classic works of Jindani et al. [40] showed  that a bactericidal treatment reduced the number of bacilli 20 times during the first two days and about 200 times during the 12 days. After two weeks of treatment, the sputum of a patient contain on average 1000 times less bacilli before treatment, a number generally too low to be detected on direct examination. Latently infected individuals who received successful chemoprophylaxis can recover at a constant rate $\eta$ and enter the recovered class $R$. Those who did not received effective chemoprophylaxis progress to active TB at a  rate $\alpha(1-\theta)$  where $\alpha$ is the rate at which latently infected individuals become infectious (this value is connected
with the average time of incubations). Once in active stage of the disease, due to their own immunity, an individual may recover naturally  and will move in the class of latently infected at rate $\gamma$.  Also, after a therapy of treatment,  infectious can spontaneously recover from the disease and will enter the   recovery class $R$ at rate $\delta$. As suggested by Styblo [41], recovered  individuals can only have a  partial immunity. Hence, they can relapse from the disease with a constant rate $\xi$ and enter the  infectious class $I$, in the same time some of them can be re-infected at a rate $(1-\xi)\beta \dfrac{\rho_{I,k}}{\rho_{k}}$ and enter class $E$. The rate  for non-disease related death is $\mu$, thus, $1/\mu$ is the average lifetime.  Infectious  have addition death rate due to  disease with a rate   $d$.\\ 
This model description is summarized in the flow diagram  given below. 

According to the derivation in [24, 25] of the continuous-time formulation for the progress
of diseases on metapopulations and assuming non-limited or frequency-dependent transmission, the equations governing the dynamics of TB propagation are
\begin{equation}
\label{network1}
\left\{\begin{array}{lcl}
\dot \rho_{S,k}&=&\Lambda- \beta\displaystyle\frac{\rho_{I,k}\rho_{S,k}}{\rho_k}-\mu \rho_{S,k}-D_S\rho_{S,k}
+kD_S\sum\limits_{k'}P(k'|k)\displaystyle\frac{\rho_{S,k'}}{k'},\\
\\
\dot \rho_{E,k}&=& \beta(1-q) \displaystyle\frac{\rho_{I,k}\rho_{S,k}}{\rho_k} +  \beta\, (1-\xi)\dfrac{\rho_{I,k}\,\rho_{R,k}}{\rho_{k}} +  \gamma\rho_{I,k}-[\mu + \eta +\alpha(1-\theta)]\rho_{E,k} \\
\\
&  & - \; D_E\rho_{E,k} + kD_E\sum\limits_{k'}P(k'|k)\displaystyle\frac{\rho_{E,k'}}{k'},\\
\\
\dot \rho_{I,k}&=&   \beta q \displaystyle\frac{\rho_{I,k}\rho_{S,k}}{\rho_k}+\alpha(1-\theta)\rho_{E,k}-(\mu+d+\gamma+ \delta)\rho_{I,k} + \xi \,\rho_{R,k} \\
\\
& & - \; D_I\rho_{I,k} +kD_I\sum\limits_{k'}P(k'|k)\displaystyle\frac{\rho_{ I,k'}}{k'},\\
	\\
 \dot \rho_{R,k} &=&  - \beta\, (1-\xi)\dfrac{\rho_{E,k}\,\rho_{R,k}}{\rho_{k}} + \eta\,\rho_{E,k} + \delta \,\rho_{I,k} - (\mu +  \xi)\,\rho_{R,k} - D_R\,\rho_{R,k}
	+kD_R\sum\limits_{k'}P(k'|k)\displaystyle\frac{\rho_{R,k'}}{k'} ,	
\end{array}\right.
\end{equation}
 where  $k$ is the degree of the patches where local population live ($k=k_1,\ldots,k_{max}$), and $P(k'|k)$ is the conditional probability
that a patch of degree $k$ has a connection to a patch of degree $k'$. As in classical
reaction-diffusion processes, Eq. (\ref{network1}) expresses the time variation
of susceptible, latently infected individuals, recovered individuals and infectious as the sum of two independent contributions: reaction and diffusion. In particular, the diffusion term includes the outflow of individuals (diffusing particles) from patches of degree $k$ and the inflow of migratory individuals from the nearest patches of degree $k'$.
 For the sake of brevity, in the sequel we  consider strictly positive diffusion rates ($D_s,D_E,D_I, D_R >0$).

For limited or frequency-dependent transmission model, we simple replace in  Eq. (\ref{network1}) the transmission term $ \beta\,\displaystyle\frac{\rho_{I,k}\rho_{S,k}}{\rho_k}$ by $\beta\,\rho_{I,k}\rho_{S,k} $.

\subsection{Positively-invariant set}
\noindent

Notice that, since  births and deaths are considered in   model (\ref{network1}), the total number of individuals is not constant at the metapopulation level. More precisely, multiplying equations in system (\ref{network1}) by $p(k)$, and summing over all $k$, we have the following differential equations for $\rho_S$, $\rho_E$, $\rho_I$ and  $\rho_R$, the average number of susceptible, latently infected, infectious   and  recovered individuals  per path at time $t$, respectively,
\begin{equation}
\label{prop1}
\left\{\begin{array}{lcl}
\dot \rho_{S}&=&\Lambda- \beta\d\sum\limits_kp(k)\d\frac{\rho_{I,k}\rho_{S,k}}{\rho_k}-\mu \rho_{S}-D_S\rho_{S}
+D_S\d\sum\limits_k\sum\limits_{k'}kp(k)P(k'|k)\d\frac{\rho_{S,k'}}{k'},\\
\\
\dot \rho_{E}&=&\beta(1-q) \d\sum\limits_kp(k)  \d\frac{\rho_{I,k}\rho_{S,k}}{\rho_k} + \beta(1-\xi) \d\sum\limits_kp(k)  \d\frac{\rho_{I,k}\rho_{R,k}}{\rho_k}  + \gamma\rho_{I}-[\mu + \eta +\alpha(1-\theta)]\rho_{E}
  \\
\\
&&- \; D_E\rho_{E} + D_E \d\sum\limits_k\sum\limits_{k'}kp(k)P(k'|k)\d\frac{\rho_{E,k'}}{k'} ,\\
\\
\dot \rho_{I}&=& \beta q \d\sum\limits_kp(k) \d\frac{\rho_{I,k}\rho_{S,k}}{\rho_k}+ \alpha(1-\theta)\rho_{E}-
(\mu+d+\gamma+ \delta)\rho_{I}\\
\\
& +&  - D_I\rho_{I}
+D_I\d \sum\limits_k\sum\limits_{k'}kp(k)P(k'|k)\d\frac{\rho_{ I,k'}}{k'} + \xi\,\rho_R ,\\
\\
  \dot \rho_{R} &=& - \beta(1-\xi) \d\sum\limits_kp(k)  \d\frac{\rho_{I,k}\rho_{R,k}}{\rho_k}  + \eta \,\rho_E +  \delta \,\rho_{I}-
(\mu +  \xi)\,\rho_{R} \\
\\
&& -\; D_R\rho_{R}
+D_R\d\sum\limits_k\sum\limits_{k'}kp(k)P(k'|k)\d\frac{\rho_{R,k'}}{k'},
\end{array}\right.
\end{equation}
where $\rho_j(t)=\d\sum\limits_kp(k)\rho_{j,k}$, $j=S,E,I,R$. Now, since the number of links emanating from nodes of degree $k$ to nodes of degree $k'$ must be equal to the number of links emanating from nodes of degree $k'$ to nodes of degree $k$ in non-directed graphs, we have the following relationship between $p(k)$ and $P(k'|k)$ [27]:
\begin{equation}
 \label{prop2}
kP(k'|k)p(k)=k'P(k|k')p(k').
\end{equation}
Using this restriction and the fact that $\sum\limits_{k}P(k|k')=1$ after changing the order of summations, Eq. (\ref{prop1}) becomes
\begin{equation}
\label{prop3}
\left\{\begin{array}{lcl}
\dot \rho_{S}&=&\Lambda- \beta\sum\limits_kp(k)\displaystyle\frac{\rho_{I,k}\rho_{S,k}}{\rho_k}-\mu \rho_{S},\\
\\
\dot \rho_{E}&=&\beta(1-q) \sum\limits_kp(k)\displaystyle\frac{\rho_{I,k}\rho_{S,k}}{\rho_k} + \beta(1-\xi) \d\sum\limits_kp(k)  \d\frac{\rho_{I,k}\rho_{R,k}}{\rho_k} +\gamma\rho_{I}-[\mu + \eta +\alpha(1-\theta)]\rho_{E},\\
\\
\dot \rho_{I}&=& \beta q \sum\limits_kp(k)\displaystyle\frac{\rho_{I,k}\rho_{S,k}}{\rho_k}  +\alpha(1-\theta)\rho_{E}-
(\mu+d+\gamma+ \delta )\rho_{I} +  \xi \,\rho_{R}, \\
\\
 \dot \rho_{R} &=& - \beta(1-\xi) \d\sum\limits_kp(k)  \d\frac{\rho_{I,k}\rho_{R,k}}{\rho_k} + \eta \,\rho_E + \delta \,\rho_{I} - (\mu + \xi)\,\rho_{R}.
\end{array}\right.
\end{equation}
Adding the expressions in the right-hand side of the equations in system (\ref{prop3}) yields
\begin{equation}
 \label{prop4}
\displaystyle\frac{d\rho}{dt}=\Lambda-\mu\rho-d\rho_I.
\end{equation}
From  the above equation, one can deduce that $ \displaystyle\frac{d\rho}{dt}\leq \Lambda-\mu\rho$. Thus, $ \displaystyle\frac{d\rho}{dt}<0$ if $\rho>\displaystyle\frac{\Lambda}{\mu}$. Since $\displaystyle\frac{d\rho}{dt}\leq \Lambda-\mu\rho$, it can be shown that using a standard comparison theorem [28], that
$$\rho(t)\leq \rho(0)e^{-\mu t}+\displaystyle\frac{\Lambda}{\mu}(1-e^{-\mu t}).$$
If $\rho(0)\leq \displaystyle\frac{\Lambda}{\mu}$, then $\rho(t)\leq \displaystyle\frac{\Lambda}{\mu}$.

\noindent
Hence,  all feasible solutions of components
of  system    (\ref{prop3}) enters the region:

\begin{equation}
\label{ens}
\Omega=\left\{(\rho_S,\rho_E,\rho_I,  \rho_R )\in\R^{4}_{\geq 0},\,\,\, \rho(t)\leq \displaystyle\frac{\Lambda}{\mu}\right\}.
\end{equation}
Thus, it follows from Eq.~(\ref{ens}) that all possible solutions of   system
 (\ref{prop3}) will enter the region $\Omega$.  Hence, the region $\Omega$, of biological interest, is positively-invariant under the flow
induced by  system (\ref{prop3}). Further, it can be shown using the theory of permanence
[26] that all solutions on the boundary of $\Omega$ eventually enter the interior of $\Omega$.
Furthermore, in $\Omega$, the usual existence, uniqueness and continuation results hold for
  system  (\ref{prop3}). Hence, system (\ref{prop3}) is well posed mathematically and epidemiologically
and it is sufficient to consider the dynamics of the flow generated by  system (\ref{prop3})
in $\Omega$. The same conclusions on $\Omega$ hold for the simple mass action (or density-dependent) model.

For networks with a connectivity pattern defined by a set of conditional
probabilities $P(k'|k)$, we define the elements of the connectivity matrix $C$
as
$$C_{kk'}=\displaystyle\frac{k}{k'}P(k'|k).$$
Note that these elements are the average number of individuals that patches of degree $k$ receive from neighboring patches of degree $k'$ assuming that one individual leaves
each of these patches by choosing at random one of the $k'$ connections [9]. One should notice that, for those degrees $k$ that are not present in the network, $P(k'|k)=0$,
$\forall k'$. Hereafter in the paper, when talking about degrees, we implicitly
mean those degrees that are present in the network. Furthermore, the case with patches
having all the same connectivity is excluded from our considerations because, under the present approach, the model equations reduce to those of a single patch SEIR model.

%


\section{Uncorrelated networks}
\noindent

In order to obtain analytical results about the  TB metapopulation dynamics, we need
to be precise about the form of $P(k'|k)$. The easiest and usual assumption is to restrict ourselves to uncorrelated networks. In these networks, the degrees of the nodes at the ends of any given link are independent, that is, no degree-degree correlation between the connected nodes. In this case, we have that
$P(k'|k)=k'p(k')/\langle k\rangle$ which corresponds to the degree distribution
of nodes (patches) that  arrive at by following a randomly chosen link [29].
\subsection{Analysis of standard incidence (or frequency-dependent) model}
\noindent

After replacing the expression of $P(k'|k)$ into   Eq.(\ref{network1}), one obtains
the following equations for TB spread in metapopulations described by uncorrelated
networks and limited transmission:
\begin{equation}\label{network2}
\left\{\begin{array}{lcl}
\dot \rho_{S,k}&=&\Lambda- \beta\displaystyle\frac{\rho_{I,k}\rho_{S,k}}{\rho_k}-\mu \rho_{S,k}-D_S\left(\rho_{S,k}-
\displaystyle\frac{k}{\langle k\rangle }\rho_S\right),\\
\\
\dot \rho_{E,k}&=& \beta(1-q)\displaystyle\frac{\rho_{I,k}\rho_{S,k}}{\rho_k}+ + \beta(1-\xi) \d\frac{\rho_{I,k}\rho_{R,k}}{\rho_k}+\gamma\rho_{I,k}-[\mu + \eta+\alpha(1-\theta)]\rho_{E,k} \\
\\
&& -D_E\left(\rho_{E,k}-
\displaystyle\frac{k}{\langle k\rangle }\rho_E\right)  ,\\
\\
\dot \rho_{I,k}&=& \beta q\displaystyle\frac{\rho_{I,k}\rho_{S,k}}{\rho_k}+\alpha(1-\theta)\rho_{E,k}-(\mu+d+\gamma+ \delta)\rho_{I,k} -
D_I\left(\rho_{I,k} - \displaystyle\frac{k}{\langle k\rangle}\rho_I\right) + \xi \,\rho_{R,k}, \\
\\
\dot \rho_{R,k} &=&  - \beta(1-\xi) \d\frac{\rho_{I,k}\rho_{R,k}}{\rho_k} + \eta \,\rho_{E,k} + \delta \,\rho_{I,k} - (\mu + \xi)\,\rho_{R,k} -  D_R\left(\rho_{R,k} - \displaystyle\frac{k}{\langle k\rangle}\rho_R\right),	
\end{array}\right.
\end{equation}
where $\langle k\rangle=\sum\limits_{k}kp(k)$ is the average network degree.

In this form, it becomes clearer that the diffusion term is simply given by the difference between the outflow of susceptible, latently infected, infectious and recovered individuals in patches of connectivity $k$, $D_s\rho_{S,k}$, $D_E\rho_{E,k}$ , $D_I\rho_{I,k}$, and $D_R\rho_{R,k}$ and the total inflow of susceptible, latently infected, infectious  and recovered individuals across all their $k$ connections, which is $k$ times the average flow of individuals across a connection in the network, $D_S\rho_S/\langle k\rangle$,
$D_E\rho_E/\langle k\rangle$ , $D_I\rho_I/\langle k\rangle$ and $D_R\rho_R/\langle k\rangle$. Note that this average flow across a connection does not depend on the degree $k$ of the considered patch because we have assumed that the architecture of the metapopulation is described
by an uncorrelated network.

In these networks, the elements of the connectivity matrix $C$ are simply
\begin{equation}
\label{conect}
 C_{kk'}=\displaystyle\frac{kp(k')}{\langle k\rangle}.
\end{equation}
Clearly, $C$ is a rank-one matrix and has the vector with components
$v_k=k$ as eigenvector of eigenvalue 1. So, if there are $n$ different
degrees in the network, then the eigenvalues of this matrix are $\lambda=0$, with algebraic multiplicity $n-1$ and $\lambda=1$ which is a simple eigenvalue. This fact will be used in the stability of equilibria of the model.
To do this, we are going to `vectorialize' system (\ref{network2}), using the following vectors of $\R^n$:

\noindent
$S=(\rho_{S,k_1},\rho_{S,k_2},\ldots,\rho_{S,k_n})^T$, $E=(\rho_{E,k_1},\rho_{E,k_2},\ldots,\rho_{E,k_n})^T$, $I=(\rho_{I,k_1},\rho_{I,k_2},\ldots,\rho_{I,k_n})^T$, \\
 $R=(\rho_{R,k_1},\rho_{R,k_2},\ldots,\rho_{R,k_n})^T$,\, $N=(\rho_{k_1},\rho_{k_2},\ldots,\rho_{k_n})^T$  and  $\mathbb I=(1,1,\ldots,1)^T$.  If $X\in\R^n$ is a vector, we denote by $\mbox{diag}(X)$ the $n\times n$ matrix whose diagonal is given by the components of $X$.  With these notations and conventions, system (\ref{network2}) becomes
\begin{equation}
 \left\{\label{network3}
\begin{array}{lcl}
\dot S&=& \Lambda\mathbb I-\beta\mbox{diag}(N)^{-1}\mbox{diag}(I)S-(\mu+D_S)S+D_SCS,\\
\\
\dot E&=&\beta(1-q)\mbox{diag}(N)^{-1}\mbox{diag}(I)S+ \beta(1-\xi)\mbox{diag}(N)^{-1}\mbox{diag}(I)R \\
\\
 && +\; \gamma I-[\mu + \eta +\alpha(1-\theta)+D_E]E+D_ECE ,\\
\\
\dot I &=& \beta q\mbox{diag}(N)^{-1}\mbox{diag}(I)S+\alpha(1-\theta)E-(\mu+d+\gamma+\delta+D_I)I + D_IC\, I + \xi \,R ,\\
\\
\dot R &=&  - \beta(1-\xi)\mbox{diag}(N)^{-1}\mbox{diag}(I)R+\eta \,E +   \delta \,I - (\mu +  \xi + D_R)\,R + D_R\,C\,R ,
\end{array}\right.
\end{equation}
where $C$  is the connectivity matrix defined as in Eq. (\ref{conect}).

We point out that in  the case where the parameters $\beta$, $q$, $\gamma$, $\mu$, $\delta$, $\theta$, $\alpha$, $\xi$, $\eta$ and $d$
are not the same for all patches, they are replaced in system (\ref{network3})
by diagonal non-negative matrices and this does not change the fundamental structure of the system.
\subsubsection{Disease-free equilibrium (DFE) for generic networks}
\noindent

The disease-free equilibrium of   model system (\ref{network1}) are the solutions $\rho^0_{S,k}$, $\rho^0_{E,k}$ and $\rho^0_{I,k}$ to the equations:
\begin{equation}
 \label{DFE}
\left\{\begin{array}{ll}
\Lambda-\beta \displaystyle\frac{\rho^0_{I,k}\rho^0_{S,k}}{\rho^0_k}-\mu \rho^0_{S,k}-D_S\rho^0_{S,k}
+kD_S\sum\limits_{k'}P(k'|k)\displaystyle\frac{\rho^0_{S,k'}}{k'}=0,\\
\\
\beta(1-q)\displaystyle\frac{\rho^0_{I,k}\rho^0_{S,k}}{\rho^0_k}+ \beta(1-\xi)\displaystyle\frac{\rho^0_{I,k}\rho^0_{R,k}}{\rho^0_k}+\gamma\rho^0_{I,k}-[\mu + \eta +\alpha(1-\theta)]\rho^0_{E,k}\\
\\
 - D_E\rho^0_{E,k}
+kD_E\sum\limits_{k'}P(k'|k)\displaystyle\frac{\rho^0_{E,k'}}{k'} =0,\\
\\
\beta q\displaystyle\frac{\rho^0_{I,k}\rho^0_{S,k}}{\rho^0_k}+ \alpha(1-\theta)\rho^0_{E,k}-(\mu+d+\gamma+\delta)\rho^0_{I,k} - D_I\rho^0_{I,k}
+kD_I\sum\limits_{k'}P(k'|k)\displaystyle\frac{\rho^0_{S,k'}}{k'} + \xi \,\rho^0_{R,k} = 0,\\
\\
 -\beta(1-\xi)\displaystyle\frac{\rho^0_{I,k}\rho^0_{R,k}}{\rho^0_k}+\eta \,\rho_{E,k} + \delta \,\rho^0_{I,k} - (\mu + \xi)\,\rho^0_{R,k} - D_R\,\rho^0_{R,k}
	+kD_R\sum\limits_{k'}P(k'|k)\displaystyle\frac{\rho^0_{R,k'}}{k'} = 0 .
\end{array}\right.
\end{equation}

For the analysis of the infection's spread, the so-called
disease-free equilibrium  is particularly relevant. By definition, this is obtained by replacing $\rho_{I,k}=0$
in Eq.(\ref{network1}), leading to an explicit expression for the number of susceptible
individuals in patches with degree $k$ that can be written as
$$
(\mu+D_S)\rho^0_{S,k}=\Lambda+D_S\sum\limits_{k'}C_{kk'}\rho^0_{S,k'}.
$$
As $\sum\limits_{k'}P(k'|k)=1$, it  follows that, for any generic network,
one has
$$\rho^0_{S,k}=\displaystyle\frac{1}{\mu+D_S}
\left(\Lambda+D_S\displaystyle\frac{k}{\langle k\rangle}\rho^0_S\right).$$
Note that Eq. (\ref{prop4}) at the disease-free equilibrium yields
$$\rho^0=\rho^0_S=\displaystyle\frac{\Lambda}{\mu}.$$
Then, the disease-free equilibrium is given by
\begin{equation}
 \label{DFE1}
\rho^0_{S,k}=\displaystyle\frac{\Lambda}{\mu+D_S}
\left(1+\displaystyle\frac{D_S}{\mu}\displaystyle\frac{k}{\langle k\rangle}\right),\qquad
\rho^0_{E,k}=\rho^0_{I,k} =  \rho^0_{R,k} =0,\qquad \forall k.
\end{equation}

Equation~(\ref{DFE1}) is also the disease free equilibrium for the simple mass action  transmission  model.
\subsubsection{Basic reproduction number and local stability of (DFE)}
\noindent

The global behavior for this model crucially depends on the basic reproduction number, that is, the average number of secondary
cases produced by a single infective individual which is introduced into an entirely susceptible population.  System~(\ref{network2}) has an evident equilibrium
$Q_0=(S^0,0,0,0)$ with $S^0_k = \rho^0_{S,k}$ defined as in Eq.~(\ref{DFE1}) and $0$ the zero vector of dimension $n$ when there is no disease. We calculate the basic
reproduction number, $\mathcal R_0$, using the next generation
approach, developed in  Ref.~[30].

Using the notations in Ref.~[30], the matrices $F$ and $V$, for the new infections and the remaining transfers, are, respectively, given by

$$F= \begin{bmatrix}
0& F_1&0\\
0& F_2&0\\
0&0&0
\end{bmatrix}\qquad \mbox{and}\qquad V= \begin{bmatrix}
 A_E\, I_n -  D_E \,C &  - \gamma I_n & 0\\
 - \alpha\,(1 - \theta)\,I_n &  A_I\, I_n -  D_I \, C & -\xi I_n \\
-\eta I_n & - \delta I_n & A_R\,I_n - D_R\,C
\end{bmatrix},$$
where $I_n$ is the identity matrix of dimension $n$,
$$F_1=\beta(1-q)I_n,\:\: F_2=\beta qI_n,\:\: A_E=[\mu + \eta + \alpha\,(1 - \theta) +  D_E ], \:\:A_I=\mu + d + \gamma+\delta+ D_I, \:\:  A_R= \mu +\xi + D_R. $$
Set
$$ F= \begin{bmatrix}
 F_{11} &  F_{12}\\
 F_{21}& F_{22}
\end{bmatrix},$$
 where
 $$ F_{11}= \begin{bmatrix}
0 &  F_1\\
0&  F_2
\end{bmatrix},\qquad  F_{12}= \begin{bmatrix}
 0\\
0
\end{bmatrix},\qquad    F_{21}= \begin{bmatrix}
 0 ,&0
\end{bmatrix}\qquad \mbox{and}\qquad F_{22}= 0.$$
Also, let
$$ V= \begin{bmatrix}
 V_{1} &  V_{2}\\
 V_{3}& V_{4}
\end{bmatrix}, $$
where
 $$V_1= \begin{bmatrix}
  A_E\, I_n -  D_E \,C &  - \gamma I_n\\
 - \alpha\,(1 - \theta)\,I_n &  A_I\, I_n -  D_I \, C
\end{bmatrix},\quad  V_2= \begin{bmatrix}
0\\
 -\xi I_n
\end{bmatrix},\quad V_3= \begin{bmatrix}
 -\eta\,I_n ,&-\delta\, I_n
\end{bmatrix}\quad \mbox{and}\quad V_4= A_R\,I_n - D_R\,C .$$
We stress that since  $V$ is a M-matrix and $-V$ is stable , then  one can deduce that $V^{-1} \geq 0 $. 

Now, we need to compute the inverse of the matrix $V$. To this end,  suppose that the inverse matrix of $V$ can be written in the following form:
$$ V^{-1}= \begin{bmatrix}
 W_{11} &  W_{12}\\
 W_{21}& W_{22}
\end{bmatrix},$$
where $W_{11}$ and $W_{22}$ are square  matrices of dimension $(2n\times 2n)$ and  $(n\times n)$, respectively.

\noindent
Observe that
$$FV^{-1}=  \begin{bmatrix}
  A &  B\\
0& 0
\end{bmatrix},$$
where $A= F_{11}\,W_{11}$ and  $B=  F_{11}\,W_{12}$. Thus, the basic reproduction ratio is defined, following [30], as the spectral radius of the next
generation matrix, $FV^{-1}$:
\begin{equation}
 \label{taux1}
\mathcal R_0 =  \rho(FV^{-1})= \rho(A)=\rho\left(F_{11}\,W_{11}\right).
\end{equation}
To compute the explicit expression of the basic reproduction number, we need to compute the inverse matrix of  $V$. To this end, we need
the following lemma stated above and proved in Appendix A.
\begin{lem}\label{inv}:
Let $N$ be a square block matrix of the following form:
$$  N= \begin{bmatrix}
N_1 &  N_2\\
N_3 & N_4
 \end{bmatrix},$$
where $N_1$ and $N_4$ are   square matrices.

\noindent
If $N_1$ and $D= N_4-N_3N_1^{-1}N_2$ are invertible, then the inverse matrix of $N$ is given by
 $$ N^{-1} = \begin{bmatrix}
 N_1^{-1}+N_1^{-1}N_2D^{-1}N_3N_1^{-1} &  -N_1^{-1}N_2D^{-1}\\
 \\
-D^{-1}N_3N_1^{-1}& D^{-1}
 \end{bmatrix}. $$
 \end{lem}

Note that the  matrix $V_1$ has the form of the matrix $N$ defined in Lemma \ref{inv}  with $N_1= A_E\,I_n - D_E\,C$, $N_2= - \gamma\, I_n$, $N_3 = - \alpha(1-\theta)\,I_n$ and $N_4=A_I\,I_n - D_I\,C$.

\noindent
Note also that the  matrix $V$ has the form of the matrix $N$ defined in Lemma \ref{inv}  with $N_1= V_1$,  $N_2= V_2$, $N_3 = V_3$ and $N_4= V_4$. So, if all the hypotheses in Lemma \ref{inv} are satisfied for the matrices  $V_1$ and $V_4$,  then Lemma \ref{inv}  can be used twice to compute  $V^{-1}$.

Thus using Lemma \ref{inv}, one can prove that  $V_1^{-1}$ has the following  form:
 $$V_1^{-1} = \begin{bmatrix}
  V_{11} & V_{12}\\
V_{21} & V_{22}
\end{bmatrix},$$
where
$$ \begin{array}{lcl}
 V_{11} &= &\left(A_E\,I_n - D_E\,C\right)^{-1} + \gamma\,\left(A_E\,I_n - D_E\,C\right)^{-1}\,\,V_{21}, \\
\\
V_{12} &= &  \gamma \left(A_E\,I_n - D_E\,C\right)^{-1}\,\,V_{22},\\
\\
 V_{21} &=  &\alpha (1- \theta) \,V_{22}\left(A_E\,I_n - D_E\,C\right)^{-1},\\
 \\
 V_{22} &=  &\left[A_I\,I_n - D_I\,C - \alpha\,\gamma(1- \theta)\left(A_E\,I_n - D_E\,C\right)^{-1}\right]^{-1}.
\end{array}$$
From the above expressions, it appears that to compute the explicit expressions of   $V_{11}$,  $ V_{12}$, $V_{21}$ and $V_{22}$,  we need to compute the inverse  matrices of 
$ \left[A_I\,I_n - D_I\,C - \alpha\,\gamma(1- \theta)\left(A_E\,I_n - D_E\,C\right)^{-1}\right]$ and  $\left(A_E\,I_n - D_E\,C\right)$. To do so, we shall used the following Lemma \ref{sam00}  stated below and proved in Appendix C.

\begin{lem}\label{sam00}: Let  $G=U+X\,W\,Z $ be an $n\times n$ invertible matrix. Suppose that the matrices $U$, $W$ and $W^{-1}+Z\,U^{-1}\,X$ are invertible. Then,
the inverse matrix of $R$ is defined as
\begin{equation}
 \label{inverse}
G^{-1}= U^{-1} -  U^{-1}\,X\,[W^{-1}+Z\,U^{-1}\,X]^{-1}\,Z\,U^{-1}.
\end{equation}
 \end{lem}

Using the above Lemma \ref{sam00} and  the fact that $C^m=C, \forall m\in \mathbb N^*$, one can easily prove that
$$ \begin{array}{ll}
\left(A_E\,I_n - D_E\,C\right)^{-1}=\displaystyle\frac{1}{A_E}\, \left[ I_n+ \displaystyle\frac{\,D_E}{A_E - D_E}\,C\right], \\
\\
\left[A_I\,I_n - D_I\,C - \alpha\,\gamma(1- \theta)\left(A_E\,I_n - D_E\,C\right)^{-1}\right]^{-1}=\displaystyle\frac{1}{a}\, \left[ I_n+ \displaystyle\frac{b}{a - b}\,C\right],
\end{array}
$$
where
$$\begin{array}{lcl}
 a= \displaystyle\frac{A_I \,(\mu+\eta+D_E)+ \alpha\,(1-\theta)(\mu+d +\delta+D_I)}{A_E}\qquad
\mbox{and} \qquad b=\displaystyle\frac{A_ED_I[\mu+\eta+\alpha(1-\theta)]+\gamma\,\alpha (1-\theta)\,D_E}{A_E[\mu+\eta+\alpha(1-\theta)]}.
\end{array}$$%
With this in mind,  after some substitutions, one has:
$$
\begin{array}{lcl}
V_{22} &= & \displaystyle\frac{1}{a}\, \left[ I_n+ \displaystyle\frac{b}{a - b}\,C\right] \\
\\
&=&  a_0\,I_n + b_0\,C , \\
\\
V_{21} &= & \displaystyle\frac{ \alpha(1-\theta)}{a\,A_E}\left[ I_n + \displaystyle\frac{b[\mu+\alpha(1-\theta)]+ a D_E}{(a-b)[\mu+\alpha(1-\theta)]}C\right] \\
 \\
 &=& a_1\,I_n + b_1\,C , \\
\\
V_{12} &= & \displaystyle\frac{ \gamma}{a\,A_E}\,\left[ I_n + \displaystyle\frac{b[\mu+\alpha(1-\theta)]+ a D_E}{(a-b)[\mu+\alpha(1-\theta)]}C\right], \\
\\
&=& a_2\,I_n + b_2\,C , \\
\\
V_{11} &= & \displaystyle\frac{a\,A_E + \gamma\,\alpha(1-\theta)}{a\,A_E^2}\left[ I_n + \displaystyle\frac{\gamma\,\alpha\,(1-\theta)\,A_E\,[b\,[\mu+\alpha(1-\theta)]+ a D_E] }{(a-b)\, [a\,A_E + \gamma\,\alpha(1-\theta)]\,[\mu+\alpha(1-\theta)]^2} C\right]\\
\\
 &+& \displaystyle\frac{a\,A_E + \gamma\,\alpha(1-\theta)}{a\,A_E^2}\left[ \displaystyle\frac{D_E\,[\mu+\alpha(1-\theta)]\, [ a\,A_E + \gamma\,\alpha(1-\theta)] }{[a\,A_E + \gamma\,\alpha(1-\theta)]\,[\mu+\alpha(1-\theta)]^2}C\right]\\
 \\
 &=&  a_3\,I_n + b_3\,C ,
\end{array}$$
where, 
$$
\begin{array}{lcl}
a_0 &=&  \dfrac{1}{a} = \d\frac{A_E}{A_I \,(\mu+\eta+D_E)+ \alpha\,(1-\theta)(\mu+d +\delta+D_I)} , \\
\\
b_0 &= & \displaystyle\frac{b}{a\,(a - b)} ,\\
 \\
 a_1&=& \displaystyle\frac{ \alpha(1-\theta)}{a\,A_E} , \\
\\
b_1 &= &  \displaystyle\frac{ \alpha(1-\theta)}{a\,A_E}\, \displaystyle\frac{b[\mu+\alpha(1-\theta)]+ a D_E}{(a-b)[\mu+\alpha(1-\theta)]}, \\
\\
a_2 & = & \displaystyle\frac{1}{a\,A_E}, \\
\\
b_2 &= & \displaystyle\frac{1}{a\,A_E}\, \displaystyle\frac{b[\mu+\alpha(1-\theta)]+ a D_E}{(a-b)[\mu+\alpha(1-\theta)]} ,\\
\\
 a_3 &=& \displaystyle\frac{a\,A_E + \gamma\,\alpha(1-\theta)}{a\,A_E^2},\\
 \\
 b_3 &=& \displaystyle\frac{a\,A_E + \gamma\,\alpha(1-\theta)}{a\,A_E^2} \left[ \displaystyle\frac{\gamma\,\alpha\,(1-\theta)\,A_E\,[b\,[\mu+\alpha(1-\theta)]+ a D_E] + (a-b)\,D_E\,[\mu+\alpha(1-\theta)]\, [ a\,A_E + \gamma\,\alpha(1-\theta)]}{(a-b)\, [a\,A_E + \gamma\,\alpha(1-\theta)]\,[\mu+\alpha(1-\theta)]^2}  \right].
\end{array}$$
This achieve the computation of $V_1^{-1}$. 

Now, we need to compute  $V^{-1}$. To this end, we need to prove the invertibility of matrix $D= V_4- V_3\,V_1^{-1}\,V_2$. Simple substitutions show that:
$$\begin{array}{lcl} D&=& V_4 - (\eta \,\xi\,V_{11} + \delta\,\xi \,V_{21}),\\
 \\
 &=& \left[ A_R - \xi \,(\eta \,a_3 + \delta \,a_1)\right]I_n - \left[ D_R + \xi \,(\eta \,b_3 + \delta\, b_1)\right]\, C.
 \end{array}
 $$
 Applying
  Lemma \ref{sam00} one again, the inverse of $D$ is given by
  $$\begin{array}{lcl} 
  D^{-1} &= & \dfrac{1}{\left[ A_R - \xi \,(\eta \,a_3 + \delta \,a_1)\right]} \left[ I_n +\dfrac{D_R + \xi \,(\eta \,b_3 + \delta\, b_1)}{ \left[A_R - \xi \,(\eta \,a_3 + \delta \,a_1)\right] - \left[ D_R + \xi \,(\eta \,b_3 + \delta\, b_1) \right]} C \right],\\
 \\
 & =& a_4\,I_n - b_4 \,C.
 \end{array}
 $$
where
 $$
 \begin{array}{lcl}
a_4 &=& \dfrac{1}{\left[ A_R - \xi \,(\eta \,a_3 + \delta \,a_1)\right]}, \\
\\
b_4 &=&  \dfrac{1}{\left[ A_R - \xi \,(\eta \,a_3 + \delta \,a_1)\right]}\dfrac{D_R + \xi \,(\eta \,b_3 + \delta\, b_1)}{ \left[A_R - \xi \,(\eta \,a_3 + \delta \,a_1)\right] - \left[ D_R + \xi \,(\eta \,b_3 + \delta\, b_1) \right]}.
\end{array}
 $$
  \\
Since $V_1$ and $D$ are invertible matrices, applying Lemma \ref{inv}, after simple calculations we have:
$$ \begin{array}{lcl}
W_{11} &=& V_1^{-1} +  V_1^{-1}\,V_2\,D^{-1}\,V_3 \,V_1^{-1}, \\
\\
&=& \begin{bmatrix}
V_{11} + \xi\,V_{12}D^{-1}\,(\eta\,V_{11} + \delta\,V_{21}) & V_{12} + \xi\,V_{12}D^{-1}\,(\eta\,V_{12} + \delta\,V_{22})\\
  \\
V_{21} + \xi\,V_{22}D^{-1}\,( \eta\,V_{11} + \delta\,V_{21}) & V_{22} + \xi\,V_{22}D^{-1}\,( \eta\,V_{12} + \delta\,V_{22)}
\end{bmatrix}.
\end{array}$$
At this stage, we need to compute the expression of $A$. Note that $A$ can be written as follows:
$$\begin{array}{lcl}
A &= &F_{11}\,W_{11}, \\
\\
& =  &  \begin{bmatrix}
F_1\left[V_{21} + \xi\,V_{22}D^{-1}\,(\eta\,V_{11} + \delta\,V_{21})\right] & F_1 \left[V_{22} + \xi\,V_{22}D^{-1}\,(\eta\,V_{12} + \delta\,V_{22}) \right] \\
  \\
F_2\left[V_{21} + \xi\,V_{22}D^{-1}\,( \eta\,V_{11} + \delta\,V_{21})\right] & F_2\left[V_{22} + \xi\,V_{22}D^{-1}\,( \eta\,V_{12} + \delta\,V_{22)}\right]
\end{bmatrix}.
\end{array}$$
%
 On the other hand, to have the explicit expression of the basic reproduction ratio, we need the following lemma whose proof is given in Appendix B.
\begin{lem} \label{det}:
Let $M$ be a square block matrix of the following form:
$$  M= \begin{bmatrix}
M_1 &  M_2\\
M_3 & M_4
 \end{bmatrix},$$
where $M_1$, $M_2 $, $M_3$  and $M_4$ are also square matrices.
\begin{enumerate}
 \item If $M_2 $ is invertible and $M_2M_3 - M_2M_4M_2^{-1}M_1=0$, then
\begin{equation}
 \label{cond1}
\rho(M) = \max \{ 0,\, \rho(M_1+M_2M_4M_2^{-1})\}.
\end{equation}
\item Moreover,  if  $M_2M_4=M_4M_2$, then
\begin{equation}
 \label{cond2}
\rho(M) = \max \{ 0,\, \rho(M_1+M_2)\}.
\end{equation}
\end{enumerate}
\end{lem}

Note that $A= F_{11}W_{11}$ has the form of the matrix $M$ defined in Lemma \ref{det} 
with \\ $M_1 = F_1\left[V_{21} + \xi\,V_{22}D^{-1}\,(\eta\,V_{11} + \delta\,V_{21})\right] $,
 $ M_2 = F_1 \left[V_{22} + \xi\,V_{22}D^{-1}\,(\eta\,V_{12} + \delta\,V_{22}) \right]$, \\Ê
 $ M_3 = F_2\left[V_{21} + \xi\,V_{22}D^{-1}\,( \eta\,V_{11} + \delta\,V_{21})\right] $  and $M_4 =  F_2\left[V_{22} + \xi\,V_{22}D^{-1}\,( \eta\,V_{12} + \delta\,V_{22)}\right]$.\\
   Since $F_1= \beta\,(1-q)\,I_n$ and $F_2= \beta\,q\,I_n$ are diagonal matrices, one has,
$$\begin{array}{lcl}
 M_2\,M_4 &=  & F_1 \left[V_{22} + \xi\,V_{22}D^{-1}\,(\eta\,V_{12} + \delta\,V_{22}) \right]\, F_2\left[V_{22} + \xi\,V_{22}D^{-1}\,( \eta\,V_{12} + \delta\,V_{22)}\right] \\
 \\
 &=& F_1\,F_2\left[V_{22} + \xi\,V_{22}D^{-1}\,(\eta\,V_{12} + \delta\,V_{22}) \right]\,\left[V_{22} + \xi\,V_{22}D^{-1}\,( \eta\,V_{12} + \delta\,V_{22)}\right] ,\\
 \\
 &=&F_2\,F_1 \,\left[V_{22} + \xi\,V_{22}D^{-1}\,( \eta\,V_{12} + \delta\,V_{22)}\right] \,\left[V_{22} + \xi\,V_{22}D^{-1}\,(\eta\,V_{12} + \delta\,V_{22}) \right] \\
\\
& = & F_2\left[V_{22} + \xi\,V_{22}D^{-1}\,(\eta\,V_{12} + \delta\,V_{22}) \right] \,F_1\left[V_{22} + \xi\,V_{22}D^{-1}\,( \eta\,V_{12} + \delta\,V_{22)}\right] ,\\
\\
 &=& M_4\,M_2,
\end{array}
$$
and
$$
 \begin{array}{lcl}
M_2\,M_3 - M_2\,M_4\,M_2^{-1}\,M_1 &= & M_2\,M_3 - M_4\,M_2\,M_2^{-1}\,M_1, \\
\\
&= &M_2\,M_3 - M_4\,M_1, \\
\\
&=  & F_1 \left[V_{22} + \xi\,V_{22}D^{-1}\,(\eta\,V_{12} + \delta\,V_{22}) \right]\,F_2\left[V_{21} + \xi\,V_{22}D^{-1}\,( \eta\,V_{11} + \delta\,V_{21})\right] - \\
\\
&&   F_2\left[V_{22} + \xi\,V_{22}D^{-1}\,( \eta\,V_{12} + \delta\,V_{22)}\right]\,F_1\left[V_{21} + \xi\,V_{22}D^{-1}\,(\eta\,V_{11} + \delta\,V_{21})\right]  , \\
\\
&=& F_1\,F_2\left[V_{22} + \xi\,V_{22}D^{-1}\,(\eta\,V_{12} + \delta\,V_{22}) \right]\,\left[V_{21} + \xi\,V_{22}D^{-1}\,( \eta\,V_{11} + \delta\,V_{21})\right] - \\
\\
 && F_2\,F_1 \left[V_{22} + \xi\,V_{22}D^{-1}\,(\eta\,V_{12} + \delta\,V_{22}) \right]\,\left[V_{21} + \xi\,V_{22}D^{-1}\,( \eta\,V_{11} + \delta\,V_{21})\right] , \\
\\
& =& 0.
\end{array}$$
With this in mind, since $A>0$, by applying Lemma \ref{det}, Eq. (\ref{taux1}) becomes
\begin{equation}
 \label{taux2}
 \begin{array}{lcl}
\mathcal R_0 &=&\rho\left[M_1 + M_2 \right],\\
 \\
 &=& \rho\left( F_1\left[V_{21} + \xi\,V_{22}D^{-1}\,(\eta\,V_{11} + \delta\,V_{21})\right] +  F_2\left[V_{22} + \xi\,V_{22}D^{-1}\,( \eta\,V_{12} + \delta\,V_{22)}\right]  \right),\\
  \\
  &=&  \beta\,\rho\left[ (1-q)\left[V_{21} + \xi\,V_{22}D^{-1}\,(\eta\,V_{11} + \delta\,V_{21})\right] +  q\left[V_{22} + \xi\,V_{22}D^{-1}\,( \eta\,V_{12} + \delta\,V_{22)}\right]  \right],\\
  \\
  &=& \beta\,\rho\left[ (1-q)V_{21}+ qV_{22} + \xi\,V_{22}D^{-1} \left [\eta\,((1-q)V_{11} + qV_{12}) +  \delta ((1-q)V_{21} + qV_{22})\right ] \right] ,\\
    \\
   &=& \beta\,\rho\left[ \left(I_n +  \delta\,\xi\,V_{22}D^{-1} \right)\left( (1-q)\,V_{21} +  q\,V_{22}\right) + \xi\,\delta\,V_{22}D^{-1}\left((1-q)V_{11} + q V_{12} \right) \right].
  \end{array}
\end{equation}
 From the above expressions, it is evident that $V_{11}, V_{12}, V_{21}, V_{22}, D^{-1}>0$.
 We point out that  as $V_{11}$ and  $V_{22}$ are irreducible  and nonnegative, one has  $V_{12}, V_{21} >0$. This implies that  $A$ is also irreducible and non-negative. Then, using the Perron-Frobenius theorem [31], one can deduce that $\rho(A)$ is a positive  eigenvalue   of $A$.  Additionally, a simple calculation can prove that
$$ \begin{array}{rcl}
\left[(1-q)\,V_{21} + q\,V_{22} \right]&= & \left[ \dfrac{(1-q)\,\alpha\,(1-\theta) + q\,A_E}{a\,A_E} \right] I_n \;+ \\
\\
&&   \left[  \dfrac{(1-q)\,\alpha\,(1-\theta)\left[ a\,D_E+ b\,(A_E- D_E)\right] + q\,b\,A_E\,(A_E-D_E)}{a\,(a-b)\,A_E\,(A_E-D_E)}\right] C\\
 \\
&= &  a_5\,I_n + b_5\,C,\\
\\
I_n + \delta\,\xi \,V_{22}\,D^{-1}  &= & I_n + \xi\,\delta\, (a_0\,I_n + b_0\,C)(a_4\,I_n + b_4\,C)\\
\\
&=&  [(1+ \xi\,\delta\,a_0a_4)\,I_n + \xi\,\delta\, (a_0b_4 + b_0a_4 + b_0b_4)C]\\
\\
&= &  a_6 \,I_n + b_6\, C , \\
\\
\xi\,\delta\,V_{22}D^{-1}\left[ (1-q)V_{11} + q V_{12}\right] &=& [(a_6 -1)\,I_n + b_6\,C] \left[ ((1-q)a_3 +qa_2 ) I_n + ( (1-q)b_3 + qb_2)\right] \\
\\
&=& (a_6 -1)\,[(1-q)a_3 +qa_2 ]I_n  \\
\\
&+& \left[ (a_6 -1)\,[(1-q)b_3 +qb_2 ] + b_6 [(1-q)a_3 +qa_2 ] + b_6 [(1-q)b_3 +qb_2] \right] C \\
\\
&=&  a_7\, I_n + b_7\, C.
\end{array}$$
Finally 
$$ \begin{array}{l}
\beta\,\left(I_n +  \delta\,\xi\,V_{22}D^{-1} \right)\left[ (1-q)\,V_{21} +  q\,V_{22}\right]+ \xi\,\delta\,V_{22}D^{-1}\left[ (1-q)V_{11} + q V_{12}\right] \\
\\
=  \beta\, [ ( a_5 \,I_n + b_5\, C)(a_6 \,I_n + b_6\, C) + \: a_7 \,I_n + b_7\, C ] \\
 \\
=  \beta\,\left[(a_5\,a_6 + a_7)\,I_n  +\left( a_5\,b_6 + b_5\,a_6 + b_5\,b_6 + b_7 \right)C \right] \\
\\
=  a_8\,I_n + b_8\,C . 
\end{array}$$
where
\begin{equation}
\label{taux-coef}
\begin{array}{lcl}
a_5 &= &\d\frac{(1-q)\,\alpha\,(1-\theta) + q\,A_E}{a\,A_E},\\
\\
b_5 &= &\d\frac{(1-q)\,\alpha\,(1-\theta)\left[ a\,D_E+ b\,(A_E- D_E)\right] + q\,b\,A_E\,(A_E-D_E)}{a\,(a-b)\,A_E\,(A_E-D_E)},\\
\\
a_6 &= & (1+ \xi\,\delta\,a_0a_4) \\
\\
b_6 &=& \xi\,\delta\, (a_0b_4 + b_0a_4 + b_0b_4), \\
\\
a_7 &=& (a_6 -1)\,[(1-q)a_3 +qa_2 ], \\
\\
b_7 &=&  \left[ (a_6 -1)\,[(1-q)b_3 +qb_2 ] + b_6 [(1-q)a_3 +qa_2 ] + b_6 [(1-q)b_3 +qb_2] \right], \\
\\
a_8 & =& a_5\,a_6 + a_7 , \\
\\
b_8 & =&  a_5\,b_6 + b_5\,a_6 + b_5\,b_6 + b_7  .
\end{array}
\end{equation}
Now, since  $C$ is a rank-one matrix that admits 1 as a unique positive eigenvalue, the greatest eigenvalue of the matrix  is $  \beta\,\left[(a_8\,I_n  + a_8\,C \right]$ is 
$  \beta\,\left[  a_8 + b_8  \right] $ and consequently, the basic reproduction ratio  of system (\ref{network2}) is
\begin{equation}
 \label{taux3}
\begin{array}{lcl}
\mathcal R_0 & =&  \beta\,\left[  a_8 + b_8  \right], \\
\\
& = & \beta\,[(a_5 + b_5)(a_6 + b_6) + a_7 + b_7 \,].
\end{array}
\end{equation}
\noindent \hfill$\square$

 The following result is established from Theorem 2 of [30]:
\begin{lem}\label{DFE}: The disease-free equilibrium $Q_0$ of system  (\ref{network2}) is locally asymptotically stable
whenever  $\mathcal R_0 < 1$, and instable if $\mathcal R_0>1$.
\end{lem}

Biologically speaking, Lemma \ref{DFE} implies that TB can be eliminated from the community (when $\mathcal R_0\leq 1$) if the initial sizes of
the population are in the basin of attraction of the disease-free equilibrium $Q_0$.

Now, let us analyze the basic reproduction number (\ref{taux3}). The parameter values  used for numerical simulation are given in Table 1.

\begin{center}
{\bf Table 1}: Description  of parameters of  model system \\
\begin{tabular}{llll}
\hline \hline Parameter &Description& Estimated value&Source\\
\hline \hline
$\Lambda$ &Recruitment rate & $1001$ year$^{-1}$& [35]\\
$\beta$&Transmission coefficient &Variable &\\
$\mu$&Per capita naturally death rate & $0.017$ year$^{-1}$&[34]\\
$q$&Fast route to active TB&$0.015$&[36]\\
$\alpha$ &Slow route to active TB &$0.0024$ year$^{-1}$ &[35]\\
$\theta$&Per capital rate of effective chemoprophylaxis&$0.001$ year$^{-1}$&[36]\\
$\delta$&Recovery rate of  infectious& $0.7372$ year$^{-1}$&[35]\\
$\eta$& \r Recovery rate due to chemoprophylaxis &year$^{-1}$&[35] \\
$\gamma$& Natural recovery rate of  infectious&$0.7372/4$ year$^{-1}$&Assumed\\
$d$&Per capita disease-induced mortality rate&$0.0012$ year$^{-1}$&[35]\\
$\xi$&Relapse of recovered individuals & $0.0986$  year$^{-1}$&[35]\\
 \hline\hline
\end{tabular}
\end{center}
\noindent

Figure \ref{Fig.1} shows the effects  of the transmission rate $\beta$  and
the patch connectivity $k$ on the basic reproduction ratio $\mathcal R_0$ given as in Eq.~(\ref{taux3}). We have taken a metapopulation with scale-free distribution $p(k)\sim k^{-3}$ with $\langle k \rangle=6$,  $k_{min}=3$ and  $D_S=D_E=D_I=D_R=1$. All other parameters are as in Table 1. The  part above the unity  of the picture corresponds to the region of the instability of the disease-free equilibrium, while the part below the unity of the figure represents the region for the stability of the disease-free equilibrium. From this figure,
one can see that $\mathcal R_0$ decreases
if $\beta$ decreases  even in the case of large values of $k$.  This means that  if the transmission coefficient $\beta$ is sufficiently small, TB infection could be
eliminated in the host population even if the number of the patch connectivity $k$ is large.   However, it is difficult to control $\beta$.
This figure also shows  that  for the chosen parameter values, if  the patch connectivity $k$ does not exceed $1.2$ ($k<6$), then TB can be controlled irrespective of the value of $\beta$. The infection will equally persist for $k>6$.

%

The combined effects   of the patch connectivity $k$ and the recovery rate  $\delta$ on the basic reproduction number $\mathcal R_0$ when $\beta=0.0017$  are shown in Fig.~\ref{Fig.2}.
This figure suggests that the basic reproduction ratio $\mathcal R_0$ decreases if $\delta$ increases or $k$ decreases.
Thus, the  treatment of TB    will have beneficial effects on   infectious populations if the   recovery rate is large.

%

%
\subsection{Analysis of the simple mass action (or density-dependent) model}
\noindent

In this section, we consider the analysis of the spread of TB in metapopulation uncorrelated networks under the assumption of simple mass action (or density-dependent). Under these assumptions, system~(\ref{network2})  can be written as
\begin{equation}\label{network4}
\left\{\begin{array}{lcl}
\dot \rho_{S,k}&=&\Lambda- \beta\,\rho_{I,k}\rho_{S,k}-\mu \rho_{S,k}-D_S\left(\rho_{S,k}-
\displaystyle\frac{k}{\langle k\rangle }\rho_S\right),\\
\\
\dot \rho_{E,k}&=& \beta(1-q)\,\rho_{I,k}\rho_{S,k} +\beta\,(1-\xi) \rho_{I,k}\,\rho_{R,k}  +\gamma\rho_{I,k}-[\mu+ \eta+\alpha(1-\theta)]\rho_{E,k}
-D_E\left(\rho_{E,k}-
\displaystyle\frac{k}{\langle k\rangle }\rho_E\right),\\
\\
\dot \rho_{I,k}&=& \beta q\,\rho_{I,k}\rho_{S,k}+\alpha(1-\theta)\rho_{E,k}-(\mu+d+\gamma+ \delta)\rho_{I,k} -
D_I\left(\rho_{I,k} - \displaystyle\frac{k}{\langle k\rangle}\rho_I \right) +  \xi \,\rho_{R,k}, \\
\\
\dot \rho_{R,k} &=&   - \beta\,(1-\xi) \rho_{I,k}\,\rho_{R,k} + \eta\,\rho_{R,k}+ \delta \,\rho_{I,k} - (\mu +  \xi)\,\rho_{R,k} -  D_R\left(\rho_{R,k} - \displaystyle\frac{k}{\langle k\rangle}\rho_R\right).	 \end{array}\right.
\end{equation}
Using the same  notations as in Eq.~(\ref{network3}), system~(\ref{network4}) can be written in the following compact form:
\begin{equation}
 \left\{\label{network5}
\begin{array}{lcl}
\dot S&=& \Lambda\mathbb I-\beta\mbox{diag}(I)S-(\mu+D_S)S+D_SCS,\\
\\
\dot E&=&\beta(1-q)\mbox{diag}(I)S+ \beta(1-\xi)\mbox{diag}(I)R + \gamma I-[\mu + \eta +\alpha(1-\theta)+D_E]E+D_EC\,E,\\
\\
\dot I &=& \beta q\mbox{diag}(I)S+\alpha(1-\theta)E-(\mu+d+\gamma+\delta+D_I)I  + D_IC\, I +  \xi R,\\
\\
\dot R &=&  - \beta(1-\xi)\mbox{diag}(I)R + \eta\,R + \delta \,I - (\mu +  \xi + D_R)\,R + D_R\,C\,R,
\end{array}\right.
\end{equation}
where $S$, $E$, $I$,  $R$ and $\mbox{diag}(I)$ are defined as in Eq.~(\ref{network3}).
\subsubsection{Local stability of the DFE}
\noindent

We give the formulae of the basic reproduction number, $\mathcal R_0$, for the density-dependent model, using
again the next generation approach, developed in [30]. Then, derive bounds on $\mathcal R_0$ in term of the
connectivities of patches.

Using the notations in [30], the matrices $F$ and $V$, for the new infections and the remaining transfers, are defined analogously as for the frequency-dependent model except that
$$F_1= \beta (1- q)\mbox{diag}(S^0)\qquad \mbox{and}\qquad  F_2= \beta q \mbox{diag}(S^0).$$
Similarly, the techniques used in the previous subsection can be used to compute the basic reproduction number of the density-dependent model~(\ref{network4}). Hence,   Lemmas~\ref{inv},   \ref{sam00} and \ref{det} can be used to find  the spectral radius  of the following matrix:
\begin{equation} \label{nextmatrix}
\begin{array}{lcl}
L&=& \beta\,\mbox{diag}(S^0)\left( \left(I_n +  \delta\,\xi\,V_{22}D^{-1} \right)\left[ (1-q)\,V_{21} +  q\,V_{22}\right]+ \xi\,\delta\,V_{22}D^{-1}\left[ (1-q)V_{11} + q V_{12}\right] \right),\\
\\
&=& \beta\,\left[ a_8 \,\mbox{diag}(S^0) + b_8\, \mbox{diag}(S^0)\, C \right],
\end{array}
\end{equation}
where $a_8$ and $ b_8 $ are defined as  in Eq.~(\ref{taux-coef}).
Thus
\begin{equation} \label{taux0}
 \mathcal R_0 =  \beta\,\rho\left[(a_8\,\mbox{diag}(S^0) + b_8 \,\mbox{diag}(S^0)\, C \right].
 \end{equation}
Since the spectral radius of $L$ is very difficult to compute, we shall only give some properties and estimates of its eigenvalues owing to its specific  form.

\noindent
We point out that $L$ is a sum of a diagonal  matrix $a_8 \,\mbox{diag}(S^0)$ and a rang $1$ matrix $ b_8\,\mbox{diag}(S^0)\, C$. Moreover, the diagonal elements of $a_8 \,\mbox{diag}(S^0)$ are positive and written in the increasing order.
Thus, $L$ can be considered as a diagonal matrix perturbed by a rank-one matrix.
Now, for a general interlacing theorem of eigenvalues for perturbations of a diagonal matrix by rank-one matrices [33], the eigenvalues $\lambda_{k_1}<\lambda_{k_2}<\ldots<\lambda_{k_n}=\lambda_{k_{\max}}$ of $L$ interlace with the eigenvalues $  \beta\,a_8\,S_{k_1}^0 <  \beta\,a_8\,S_{k_2}^0 < \ldots <  \beta\,a_8\,S_{k_n}^0$ of $  \beta\,a_8\,\mbox{diag}(S^0)$ as follows
$$   \beta\, a_8\,S_{k_1}^0< \lambda_{k_1}<  \beta\,a_8\,S_{k_2}^0 < \lambda_{k_2}<\ldots <  \beta\,a_8\, S_{k_n}^0 <\lambda_{k_n} = \lambda_{k_{\max}}. $$
Then, it follows that all the eigenvalues of  $L$ are real, simple, positive and the greatest one is  $\lambda_{k_{\max}}=  \mathcal R_0 $. Thus, the following inequality for  $\mathcal R_0 $ holds:
$$\mathcal R_0  >  \beta\,a_8\,\rho^0_{S,k_{\max}}.$$
 Note that   $\rho^0_{S,k_{max}}$ is defined as
$$\rho^0_{S,k_{max}}=\displaystyle\frac{\Lambda}{\mu\langle k \rangle (\mu+D_S)}[\mu\langle k \rangle+k_{max}D_S].$$
Therefore, a sufficient condition for the DFE to be unstable is given by the following lemma:
\begin{lem} \label{unstable}
If
\begin{equation} \label{cond}
  \beta\,a_8\,\d\frac{\Lambda \,( \mu\,\langle k\rangle+D_S\,k_{max} ) }{ \mu\,\langle k\rangle\,(\mu+D_S)}  > 1,
 \end{equation}
 then the disease-free equilibrium of the density dependent model is unstable.
\end{lem}
Condition~(\ref{cond}) implies that $\mathcal R_0>1$, which is a sufficient condition for the DFE to be unstable. Rearranging  condition~(\ref{cond}) gives
\begin{equation} \label{cond2}
 \rho^0_{S,k_{\max}} > \d\frac{1}{  \beta\,a_8}.
  \end{equation}

Condition~(\ref{cond2}) simply says that, if the number of individuals inhabiting those patches
with highest  connectivity in the metapopulation, for fixed values of $\mu$, $\gamma$, $D_E$, $D_I$, $D_R$, $\delta$, $\xi$, $q$,
$\theta$, $\beta$, $d$, $\gamma$, $\eta$ and $\alpha$, a large enough  $\rho^0_{S,k_{\max}}$  guarantee the instability of the disease-free equilibrium. This implies  that the infection reaches all patches.

Now,  we prove that $\mathcal R_0$ is bounded above and below and give a sufficient condition of the instability of the DFE in term of  the average density of patches of lowest connectivities.

Observe that $L$ is nonnegative ($L\geq 0$) and  $S_k^0$ is an increasing function of the connectivity $k$. Thus
$$ \beta\,\left [\,a_8 (\min_k{S_k^0})\,I_n+ b_8 \, \mbox{diag}(S^0)\, C \right] < L < \beta\,\left [a_8 \,(\max_k{S_k^0})\,I_n+  b_8\, \mbox{diag}(S^0)\, C \right].
$$
Since $S^0_k = \rho^0_{S,k} $, one has
$$ \beta\,\left [ S^0_{k_{\min}} a_8\,I_n+ b_8 \, \mbox{diag}(S^0)\, C \right] < L < \beta\,\left [ S^0_{k_{\max}} a_8\,I_n+ b_8\, \mbox{diag}(S^0)\, C \right].
$$
Then, one can deduce that
$$ \beta\,\rho \left [ S^0_{k_{\min}} a_8\, I_n+ b_8\, \mbox{diag}(S^0)\, C \right] < \rho(L) <  \beta \,\rho \left [ S^0_{k_{\max}} a_8\,I_n+ b_8\, \mbox{diag}(S^0)\, C \right],$$
which implies
$$ \beta\,\left[ a_8\,S^0_{k_{\min}}  + b_8\, \sum_k S_k^0\, C_{kk} \right] < \rho(L) <  \beta\,\left [ a_8\,S^0_{k_{\max}} + b_8\, \sum_k S_k^0\, C_{kk} \right]. $$
We have established the following lemma which give precise bounds on $\mathcal R_0$ and then yield a sufficient condition for the instability of the DFE in term of  the average density of patches of lowest connectivities.
\begin{lem} \label{bounds}: The basic reproduction number of the density-dependent model satisfies
\begin{equation} \label{R0-bounds}
\beta\,\left[ a_8\,S^0_{k_{\min}}  + b_8\, \sum_k S_k^0\,\dfrac{kp(k)}{\langle k\rangle} \right] < \mathcal R_0 <  \beta\,\left [ a_8\,S^0_{k_{\max}} + b_8\, \sum_k S_k^0\,\d\frac{kp(k)}{\langle k\rangle} \right].
\end{equation}
\end{lem}

The proof of this lemma is straightforward since  $b_8\,\mbox{diag}(S^0)C $ is  a rank one matrix, therefore, the only non zero eigenvalue of this matrix is the sum $b_8\, \d\sum_k S_k^0\,\d\frac{kp(k)}{\langle k\rangle} $ of its diagonal entries.

From this Lemma~\ref{bounds}, we deduce a sufficient condition of the instability of the DFE in term of  the average density of patches of lowest connectivities given as:
\begin{equation} \label{cond3}
 \beta\,\left[ a_8\, S^0_{k_{\min}} + b_8\, \sum_k S_k^0\,\dfrac{kp(k)}{\langle k\rangle} \right] > 1.
 \end{equation}
Since $S^0_k = \rho^0_{S,k} $, condition (\ref{cond3}) becomes
\begin{equation} \label{cond4}
 \rho^0_{S,k_{\min}} > \d\frac{1}{a_8}\left[ \dfrac{1}{\beta} - b_8 \,\d\sum_k \rho^0_{S,k}\,\d\frac{kp(k)}{\langle k\rangle} \right]
  \end{equation}
This condition~(\ref{cond4})  says that, if the number of individuals inhabiting those patches
with lowest connectivity in the metapopulation, for fixed values of $\mu$, $\gamma$, $D_E$, $D_I$, $D_R$, $\delta$, $\xi$, $q$,
$\theta$, $\beta$, $d$, $\gamma$, $\eta$ and $\alpha$, a large enough  $\rho^0_{S,k_{\min}}$  guarantee the instability of the disease-free equilibrium. This implies  that the infection reaches all patches.

In summary, it is classically known  that if $\mathcal R_0 < 1$, then the DFE is locally stable, and if $\mathcal R_0 >1$, then it is unstable. With this classic result in mind and the bounds on $\mathcal R_0$ giving by condition~(\ref{cond}) and condition~(\ref{cond4}), we have established the following result giving sufficient conditions for the instability of the DFE.
\begin{thm}:
For the model with density-dependent  model~(\ref{network4}), \\
\noindent
if the average density of patches with highest connectivities  satisfies
\begin{equation}
\label{gra}
\rho^0_{S,k_{\max}} >  \dfrac{1}{ \beta\,a_8}
\end{equation}
or, \\
\noindent
if the average density of patches with lowest  connectivities  satisfies
\begin{equation}
 \rho^0_{S,k_{\min}} > \d\frac{1}{a_8}\left[ \dfrac{1}{\beta} - b_8 \,\d\sum_k \rho^0_{S,k}\,\d\frac{kp(k)}{\langle k\rangle} \right]\end{equation}
then, the DFE is unstable.
\end{thm}

Model of this type demonstrates clear infection threshold. In the presence of a threshold, disease eradication requires the reduction of the infection rate below a critical level  where a stable infection-free equilibrium is guaranteed. In epidemiological terminology, the infection threshold may be expressed in terms of the basic reproductive ratio $\mathcal R_0$, the average number of infections produced by a single infected individual in a
population of susceptible. From this definition, it is clear  that TB infection can spread in a population only if $\mathcal R_0>1$.

In conclusion, crossing the threshold reduces the basic reproductive ratio $\mathcal R_0$ below  unity and the infection is prevented from propagating.


\subsubsection{Endemic equilibrium}
\noindent

{\r Herein, we investigate the existence of an endemic equilibrium of  system (\ref{network5})  in the special case where there is no re-infection after recovery (i.e. no flow from the recovered class to the latently infected class due to infection), but with possible relapse from the disease}. 

 To this end, it is more convenient to write system (\ref{network5}) in a more compact form. In a more compact form,   model (\ref{network5}) may be written as follows:
\begin{equation}
\label{network6}
\left\{\begin{array}{l}
\dot x= \Lambda\mathbb I-\mbox{diag}(B\,y)\,x+[D_S\,C-(\mu+D_S)]\,x,\\
\\
\dot y=\sum\limits^n_{i=1}\langle\,e_i\mid\,B\,y\rangle\langle\,e_i\mid\,x\rangle\,\mathcal K_i -V\,y,
\end{array}\right.
\end{equation}
where $x=S\in\R^{n}_{\geq 0}$, $y=(E,I, R)^T\in\R^{3n}_{\geq 0}$, $\mathcal K_i\in\R^{3n}$  are constant vectors with\\
\\
$\mathcal K_1=(\underbrace{1-q,0,\cdots,0},\underbrace{q,0,\cdots,0}, \underbrace{0,\cdots,0})^T$, \\
$\mathcal K_2=(\underbrace{0,1-q,0,\cdots,0},\underbrace{0,q,0,\cdots,0}, \underbrace{0,\cdots,0})^T,\\
\vdots$ \\
$\mathcal K_n=(\underbrace{0,\cdots,0,1-q},\underbrace{0,\cdots,0,q}, \underbrace{0,\cdots,0})^T$, 
\\
$e_i$ is the canonical basis of $\R^n$, 
$B=[0,\beta\,I_n,  0]$ with $0$  a $n\times n$ null matrix, $\mathbb I$ is defined as in Eq. (\ref{network5})
 and $V$ is the  $3n\times 3n $ constant matrix:
$$  V= \begin{bmatrix}
 A_E\, I_n -  D_E \,C &  - \gamma I_n& 0\\
 \\
 - \alpha\,(1 - \theta)\,I_n &  A_I\, I_n -  D_I \, C & -\xi I_n\\
 \\
 -\eta & - \delta I_n & A_R\,I_n - D_R\,C
\end{bmatrix}.$$

We point out that the matrix $-V$ is a Metzler matrix, that is, a matrix with all its off-diagonal entries nonnegative [31-34].

 With  this new notations, and using  the method of [30], the basic reproduction ratio (\ref{taux0}) satisfies
 \begin{equation}\label{basic-ratio-metapop}
\mathcal R_0 = \rho\left[\sum\limits^n_{i=1}\langle\,e_i\mid\,x^0\rangle\,B\,V^{-1}\,\mathcal K_i \,e_i^T\right].
\end{equation}
where $x^0= S^0 = (\rho^0_{S,k})_k$.

Let $Q^*=(x^*,y^*)$ be the positive endemic equilibrium of system (\ref%
{network6}). Then, the positive endemic equilibrium (steady state with $y>0$) can be
obtained by setting the right hand side of equations in system (\ref%
{network6}) at zero, giving
\begin{equation}
\label{equilibrium}
\left\{\begin{array}{l}
 \Lambda\mathbb I-\mbox{diag}(B\,y^*)\,x^*+[D_S\,C-(\mu+D_S)I_n]\,x^*=0,\\
\\
\sum\limits^n_{i=1}\langle\,e_i\mid\,B\,y^*\rangle\langle\,e_i\mid\,x^*\rangle\,\mathcal K_i -V \,y^*=0.
\end{array}\right.
\end{equation}
Multiplying the second equation of (\ref{equilibrium}) by $V^{-1}$ yields
$$
y^*=\sum\limits^n_{i=1}\langle\,e_i\mid\,B\, y^*\rangle\langle\,e_i\mid\, x^*\rangle\,V^{-1}\,\mathcal K_i.
$$
Using the first equation of (\ref{equilibrium}), one has
$$x^*= [\mbox{diag}(B\,y^*)-[D_S\,C-(\mu+D_S)I_n] ]^{-1}\Lambda\mathbb I.$$
Then, one can deduce that
\begin{equation}
\label{end1}
 y^*=\sum\limits^n_{i=1}\langle\,e_i\mid\,B\, y^*\rangle\langle\,e_i\mid\,[\mbox{diag}(B\, y^*)-[D_S\,C-(\mu+D_S)I_n]]^{-1}\Lambda\mathbb I\rangle\,V^{-1}\,\mathcal K_i.
\end{equation}

Remind that at the disease-free equilibrium, one has
$$\Lambda\mathbb I=-[D_S\,C-(\mu+D_S)I_n]\,x^0\gg0.$$
Plugging the above expression in  Eq. (\ref{end1}) yields
$$ y^*=\sum\limits^n_{i=1}\langle\,e_i\mid\,B\, y^*\rangle\langle\,e_i\mid\,-[\mbox{diag}(B\, y^*)-[D_S\,C-(\mu+D_S)I_n]]^{-1}[D_S\,C-(\mu+D_S)I_n]\,x^0\rangle\,V^{-1}\,\mathcal K_i.$$
Multiplying the above equation by $B$ and setting  $z^*=B y^*$ gives
\begin{equation}\label{end2}
z^*=\sum\limits^n_{i=1}\langle\,e_i\mid\,z^*\rangle\langle\,e_i\mid\,-P^{-1}( z^*)[D_S\,C-(\mu+D_S)I_n]\,x^0\rangle\,B\,V^{-1}\,\mathcal K_i,
\end{equation}
where
$$P(z^*)=\mbox{diag}( z^*)-[D_S\,C-(\mu+D_S)I_n].$$
 We give the explicit expression of the inverse matrix of $P(z^*)$ since we will need it later.  Note that $P(z^*)$  has the form of the matrix  $R = U+XWZ$ given in  Lemma \ref{sam00} with $U=\mbox{diag}[z^*+(\mu+D_S)\mathbb I]$, $X=[k_1,k_2,\ldots,k_n]$,
$W=1$ and $Z=\displaystyle\frac{D_S}{\langle k \rangle}[p(k_1),p(k_2),\ldots,p(k_n)]$.
Then, using Lemma \ref{sam00}, a simple computation gives
\begin{equation}
 \label{mat}
 P^{-1}(z^*)=\mbox{diag}\left[\frac{1}{z^*_k +\mu+D_S}\right]\,\left[ I_n + \displaystyle\frac{D_S\,C\,\mbox{diag}\left[\displaystyle\frac{1}{z^*_k +\mu+D_S}\right]}{1- \displaystyle\frac{D_S}{\langle k \rangle}\, \sum_k \,\displaystyle\frac{kp(k)}{ z^*_k+\mu+D_S}}\right].
\end{equation}
Now, from Eq. (\ref{end2}), one has
\begin{equation}\label{end3}
\langle\,e_j\mid z^*\rangle=\sum\limits^n_{i=1}\langle\,e_i\mid\, z^*\rangle\langle\,e_i\mid\,P^{-1}( z^*)P(0)\,x^0\rangle\langle\,e_j\mid\,B\,V^{-1}\,\mathcal K_i \rangle,\,\,\, j=1,2,\ldots,n,
\end{equation}
where $P(0)=-[D_S\,C-(\mu+D_S)I_n]$.  From the above equation, one can deduce  that
\begin{equation}\label{end4}
\sum\limits^n_{j=1}\langle\,e_j\mid z^*\rangle=\sum\limits^n_{i=1}\langle\,e_i\mid\, z^*\rangle\langle\,e_i\mid\,P^{-1}( z^*)P(0)\,x^0\rangle\left\langle\,\sum\limits^n_{j=1}e_j\mid\,B\,V^{-1}\,\mathcal K_i \right\rangle.
\end{equation}
Then, to find the endemic equilibrium of system (\ref{network5}), it suffices to find solutions of  the following equation:
\begin{equation}
 \label{bow}
H(z^*)=1,
\end{equation}
where
\begin{equation}
 \label{fun}
H( z^*)=\displaystyle\frac{\sum\limits^n_{i=1}\langle\,e_i\mid\,z^*\rangle\langle\,e_i\mid\,P^{-1}(z^*)P(0)\,x^0\rangle\left\langle\,\sum\limits^n_{j=1}e_j\mid\,B\,V^{-1}\,\mathcal K_i \right\rangle}{\sum\limits^n_{j=1}\langle\,e_j\mid z^*\rangle},
\end{equation}
where $P^{-1}(z^*)$ is defined as in Eq. (\ref{mat}). Note  that $z^*$ are the intersection points between the curve of $H(z^*)$ and the line $z=1$.

From Eq. (\ref{fun}), it follows that the function $H(z^*)$ satisfies
 $$\lim\limits_{z^*\rightarrow+\infty}H(z^*)=0,$$
and
$$\lim\limits_{ z^*\rightarrow0}H( z^*)=\sum\limits^n_{i=1}\langle e_i\mid x^0\rangle\,\left\langle\,\sum\limits^n_{j=1}e_j\mid\,B\,V^{-1}\,\mathcal K_i \right\rangle.$$

We claim the following result.
\begin{lem}\label{infRo}: The inequality
$\lim\limits_{ z^*\rightarrow0}H( z^*)\geq \mathcal R_0$ holds.
\end{lem}

\noindent
{\bf Proof}: Let $A= \sum\limits^n_{i=1}\langle\,e_i\mid\,x^0\rangle\,B\,V^{-1}\,\mathcal K_i \,e_i^T $. Then, using Eq. (\ref{basic-ratio-metapop}), one has  $\mathcal R_0= \rho(A)$. Since $A$ is a nonnegative matrix, if  $ r_j =\sum\limits_i^n A_{ij}$ is the sum of the $j^{th}$ column of $A$, one has
$$ \min\limits_j \{r_j\} \leq \rho(A) \leq  \max\limits_j \{r_j\}.$$
If $e_j $ denotes the canonical basis of $\mathbb R^n$, $ \mathbb I= (e_1+e_2+ \cdots + e_n)^T$, using the fact that $e_i^T \mathbb I= 1$, $\forall i$, one has
$$\begin{array}{lcl}
r_j= e_j^T\,A\, \mathbb I &=& e_j^T\left( \sum\limits^n_{i=1}\langle\,e_i\mid\,x^0\rangle\,B\,V^{-1}\,\mathcal K_i \,e_i^T \right)\,\mathbb I, \\
\\
&=& e_j^T\left( \sum\limits^n_{i=1}\langle\,e_i\mid\,x^0\rangle\,B\,V^{-1}\,\mathcal K_i \right), \\
\\
&=& \left\langle\,e_j\mid\,\sum\limits^n_{i=1} \langle e_i\mid x^0\rangle\,B\,V^{-1}\,\mathcal K_i \right\rangle, \\
\\
&=& \sum\limits^n_{i=1} \left\langle\, e_i\mid x^0 \rangle \langle\,e_j\mid\,B\,V^{-1}\,\mathcal K_i \right\rangle.
\end{array}
$$
With this mind, one can deduce that
$$\begin{array}{lcl}
\sum\limits^n_{j=1} r_j &=&\sum\limits^n_{j=1}e_j^T\,A\, \mathbb I ,\\
\\
&=&\sum\limits^n_{i=1}\langle e_i\mid x^0\rangle\,\left\langle\,\sum\limits^n_{j=1}e_j\mid\,B\,V^{-1}\,\mathcal K_i \right\rangle,
 \\
\\
&=& \lim\limits_{ z^*\rightarrow0}H( z^*).
\end{array}$$
Then, one has that
$$ \mathcal R_0 = \rho(A) \leq \max \limits_j \{r_j \} \leq \sum\limits_j^n r_j,$$
which implies that  $ \lim\limits_{ z^*\rightarrow0}H( z^*) \geq \mathcal R_0$.  This completes the proof.

\noindent \hfill$\square$

 Note that we use the expression of  $V^{-1}$ to put emphasis on the fact that $V^{-1}\geq 0$ because  $-V$ is a Metzler matrix. Since   $ \lim\limits_{ z^*\rightarrow0}H( z^*)\geq \mathcal R_0$ and $\lim\limits_{z^*\rightarrow+\infty}H(z^*)=0$,    $H(z^*)$ is a positive function. Thus, positive solutions of Eq. (\ref{bow}) exist if and only if $\lim\limits_{ z^*\rightarrow0}H( z^*)>\mathcal R_0>1$.  From the first equation of (\ref{equilibrium}), one has $x^*= P^{-1}(z^*)\Lambda\mathbb I$.  Since $P^{-1}(z^*)$ is a positive definite matrix, one has $x^*>0$. On the other hand, since $z^*$ are the intersection points between the curve of $H(z^*)$ and the line $z=1$, one has  that $z^*>0$. Then, when $\mathcal R_0>1$, the equilibria are endemic.
This means that  there exists at least one endemic equilibrium of the model (\ref{network5}). Also, note that  $ z^* = B\,y^* $ is not a bijection (it is a onto map, but not a one to one map), one can conclude  that the TB model with simple mass action transmission could have multiple endemic equilibria. However, to know the number of endemic equilibria, we need to analyze the function $H(z^*)$. We stress that   Eq.~(\ref{bow}) is very difficult to solve analytically due to the fact that   $H$ is a highly nonlinear function. Nonetheless, one can numerically plot this curve and examine how the intersection point(s) with the line $z=1$ change
with model parameters.
We have established the following theorem for the density-dependent model (\ref{network5}).
\begin{thm}:
For the model with density-dependent  model (\ref{network4}), if the basic reproduction number $\mathcal R_0 > 1$, then there exists at least one endemic equilibrium.%
\end{thm}
%
\subsubsection{Numerical studies}
\noindent

To illustrate the various theoretical results contained in the previous section,
system~(\ref{network4})  are
simulated using the parameter value/range in Table 1. In all simulations, the initial conditions have been chosen randomly.
We have also taken a metapopulation with scale-free distribution $p(k)\sim k^{-3}$ with $\langle k \rangle=6$ and   $k_{min}=3$.

Figure~\ref{Fig.3} gives the evolution of the model~(\ref{network4}) when $\beta=0.0001$ and $D_S=D_E=D_I=D_R=1$ (so that $\mathcal R_0<1$). All other parameters are as in Fig.~\ref{Fig.1}.
Figure~\ref{Fig.3}(a) presents the prevalence curves of the model while, the time evolution of the number of infected individuals in  each patch is depicted in Fig.~\ref{Fig.3}(b).
From these figures, it clearly appears  that the disease disappears in the host population even for higher values of the patch connectivity.

%
Figure~\ref{Fig.4} gives the evolution of the model~(\ref{network4}) when $\beta=0.001$ and $D_S=D_E=D_I=D_R=1$ (so that $\mathcal R_0>1$). All other parameters are as in Fig.~\ref{Fig.1}. From this figure, one can observe that   the disease persists in the host population. In addition, one can also observe that  as the patch connectivity increases, the prevalence of the infection also increases.


Now, let us examine the influence of the migration on  the propagation of TB in the host population.

Figure~\ref{Fig.5} presents the prevalence of the infection of model~(\ref{network4})  in nodes of degree $k$ of an uncorrelated  scale-free network for different values of the migration rates. From this figure,   the role of the migration rates $D_S$, $D_E$, $D_I$ and $D_R$ is remarkable. Increasing the value of the migration rates $D_S$, $D_E$, $D_I$ and $D_R$ causes a reduction in the prevalence of the infection. This is the only case we have observed in which the infection prevalence changes non-uniformly across the metapopulation when varying the value of a parameter.

%

\section{Conclusion}
\noindent

In this paper,  we have presented a system of differential equations of reaction-diffusion type describing the TB spread in heterogeneous complex  metapopulations. The spatial configuration is given by the degree $p(k)$ and the conditional probabilities $P(k\mid k')$. For uncorrelated networks under the assumption of standard incidence  transmission, we have computed the disease-free equilibrium and the basic reproduction number. We have also shown that the disease-free equilibrium is locally asymptotically stable. Moreover, for uncorrelated networks and under assumption of simple mass action  transmission, necessary and sufficient conditions for the instability of the disease-free equilibrium for uncorrelated networks have been given in term of the highest and lowest  connectivities of patches. We have also shown that  if the basic reproduction number $\mathcal R_0>1$, then the simple mass action  model could have multiple endemic equilibria. Through numerical simulations, we found that the prevalence of the infection increases with the path  connectivity. Also, increasing the value of the migration rates cause a reduction in the prevalence of the infection.

\section*{Acknowledgments}
\noindent

Berge Tsanou   acknowledges with thanks
the support  of AUF  (Agence Universitaire de la Francophonie), Bureau Afrique Centrale. Samuel Bowong   acknowledges  the support  of the Alexander von Humboldt Foundation, Germany. We are grateful to the reviewer for insightful comments.

\appendix
\section*{Appendix A: Proof of Lemma~\ref{inv} }
\noindent

In this appendix, we give the proof of Lemma~\ref{inv}. Note that the matrix $N$ can be written as
$$ \begin{array}{lcl}
 N &= & \begin{bmatrix}
N_1 &  N_2\\
\\
N_3 & N_4
 \end{bmatrix},\\
\\
&=&  \begin{bmatrix}
N_1 &  0\\
\\
N_3 & I
 \end{bmatrix} \begin{bmatrix}
I &  N_1^{-1}N_2\\
\\
0 & D
 \end{bmatrix}.
\end{array}$$
Then, one can deduce that
$$
 \begin{array}{lcl}
 N^{-1} &=  &\begin{bmatrix}
I &  N_1^{-1}N_2\\
\\
0 & D
 \end{bmatrix}^{-1}  \, \begin{bmatrix}
N_1 &  0\\
\\
N_3 & I
 \end{bmatrix}^{-1},\\
 \\
 &= & \begin{bmatrix}
I &  -N^{-1}_1N_2\,D^{-1}\\
\\
0 & D^{-1}
 \end{bmatrix}\,  \begin{bmatrix}
N_1^{-1}&  0\\
\\
-N_3N_1^{-1} & I
 \end{bmatrix},\\
 \\
 &=&\begin{bmatrix}
 N_1^{-1}+N_1^{-1}N_2D^{-1}N_3N_1^{-1} &  -N_1^{-1}N_2D^{-1}\\
 \\
-D^{-1}N_3N_1^{-1}& D^{-1}
 \end{bmatrix}.
 \end{array}$$
This ends the proof.

\noindent \hfill$\square$

\section*{Appendix B: Proof of Lemma~\ref{det} }
\noindent

In this appendix, we give the proof of Lemma~\ref{det}. To do so, we shall use the properties of the determinant.

\noindent
Let $\lambda $  the spectrum of $M$. Assume that $M$ is a $2n\times 2n$ square matrix, then,
$$ \begin{array}{lcl}
 \det(M- \lambda\,I_{2n}) &= & \det \begin{bmatrix}
M_{1}-  \lambda\,I_{n}& M_{2} \\
\\
M_{3} & M_{4}- \lambda\,I_{n}
\end{bmatrix}, \\
\\
&= & (-1)^n\,\det \begin{bmatrix}
 M_{2} & M_{1}-  \lambda\,I_{n} \\
 \\
 M_{4}- \lambda\,I_{n} & M_{3}
 \end{bmatrix},\\
 \\
 &= &(-1)^n\,\det( M_{2} )\,det\left[M_{3}-(M_{4}-\lambda I_n)\,M_{2}^{-1}\,(M_{1}-\lambda I_n)  \right], \\
 \\
 &= & (-1)^n\,\det\left[M_{2} \,M_{3}\,-M_{2} \,M_{4}\,M_{2}^{-1}\,M_{1}+\lambda \, (M_{1}+M_{2}\,M_{4}\,M_{2}^{-1}) - \lambda^2 I_n \right].
 \end{array}
$$
If $ M_{2} \,M_{3}\,-M_{2} \,M_{4}\,M_{2}^{-1}\,M_{1}=0$,  then
$$ \det(M- \lambda\,I_{2n}) = (-\lambda)^n\,\det\left[M_{1}+M_{2}\,M_{4}\,M_{2}^{-1}-
 \lambda \,I_n\right].$$
Moreover if $M_2\,M_4 = M_4\,M_2$ then
$$ \det(M- \lambda\,I_{2n}) = (-\lambda)^n\,\det\left[M_1 + M_4 - \lambda \,I_n\right].$$
This ends the proof.

\noindent \hfill$\square$
\section*{Appendix C: Proof of Lemma~\ref{sam00} }
\noindent

In this Appendix, we give the proof of Lemma\ref{sam00}. To do so, it suffices to verified that  $GG^{-1}= I_n$. Indeed,  one has
$$\begin{array}{lcl}
 GG^{-1}& =& UU^{-1} - X\left[W^{-1} + ZU^{-1}X\right]^{-1}ZU^{-1} + XWZU^{-1}\\
\\
 &-& XWZU^{-1}X\left[W^{-1}+ZU^{-1}X\right]^{-1}ZU^{-1},\\
\\
& =& I_n - X\left[\left[W^{-1} + ZU^{-1}X\right]^{-1}+W-WZU^{-1}X\left[W^{-1} + ZU^{-1}X\right]^{-1}\right]ZU^{-1},\\
\\
&=& I_n - XW\left[ W^{-1} \left[W^{-1} + ZU^{-1}X\right]^{-1} - I_n + ZU^{-1}X\left[W^{-1} + ZU^{-1}X\right]^{-1}\right]ZU^{-1},\\
 \\
& =& I_n - XW\left[  \left[W^{-1} + ZU^{-1}X\right]\left[W^{-1} + ZU^{-1}X\right]^{-1} - I_n \right]ZU^{-1},\\
\\
& =& I_n -XW(I_n-I_n)ZU^{-1},\\
 \\
& =& I_n.
\end{array}$$

This concludes the proof.

\noindent \hfill$\square$
\end{document}